\documentstyle[amssymb]{amsart}

\numberwithin{equation}{section}
\theoremstyle{plain}
\newtheorem{thm}{Theorem}[section]
\newtheorem{prop}[thm]{Proposition}

\newtheorem{conj}[thm]{Conjecture}
\newtheorem{prob}[thm]{Problem}
\theoremstyle{definition}

\theoremstyle{remark}
\newtheorem{rem}[thm]{Remark}

\newcommand{\comment}[1]{}
\newcommand{\br}[1]{\langle #1 \rangle}
\newcommand{\contfrac}[2]{\dfrac{\hfill{#1}\hfill|}{|\hfill{#2}\hfill}}
\newcommand{\BZ}{{\Bbb Z}}

\newcommand{\BR}{{\Bbb R}}
\newcommand{\BC}{{\Bbb C}}

\title[ \ \ ]
{Affine Weyl groups, 
discrete dynamical systems and 
Painlev\'e equations}
\begin{document}
\maketitle
\begin{center}
Masatoshi Noumi and Yasuhiko Yamada
\end{center}
{\small
\begin{center}
Department of Mathematics, Kobe University
Rokko, Kobe 657-8501, Japan
\end{center}
}

\begin{abstract}
A new class of representations of affine Weyl 
groups on rational functions are constructed, 
in order to formulate 
discrete dynamical systems associated 
with affine root systems. 
As an application, 
some examples of difference and differential 
systems of Painlev\'e type are discussed. 
\end{abstract}

\section*{Introduction}

In this paper, we propose a class of discrete dynamical systems 
associated with affine root systems,  by constructing new 
representations of affine Weyl groups. 
This class of difference systems covers certain types of  
discrete Painlev\'e equations, 
and is expected also to provide 
a general framework to describe the structure of B\"acklund 
transformations of differential systems of Painlev\'e type.  

\medskip
By a series of works by K.~Okamoto \cite{O}, 
it has been known since 80's 
that Painlev\'e equations 
$P_{\,\text{II}}$, $P_{\,\text{III}}$, $P_{\,\text{IV}}$, 
$P_{\,\text{V}}$ and $P_{\,\text{VI}}$ 
admit the affine Weyl groups of type 
$A^{(1)}_{1}$, $C^{(1)}_{2}$, $A^{(1)}_{2}$, $A^{(1)}_{3}$ and $D^{(1)}_{4}$, 
respectively, as groups of B\"acklund transformations. 
The relationship between the affine Weyl group symmetry 
and the structure of classical solutions has been clarified 
through the studies of irreducibility of Painlev\'e equations 
in the modern sense of H.~Umemura 
(see \cite{O},\cite{Mu},\cite{U0},\cite{NO}, for instance). 
 
In a recent work \cite{NY1}, the authors introduced a new representation 
(\ref{SP4})
of the fourth Painlev\'e equation $P_{\,\text{IV}}$ from which the structures  
of B\"acklund transformations and of special solutions of $P_{\,\text{IV}}$ 
are understood naturally.  
This sort of ``symmetric forms'' can be formulated for other Painlev\'e 
equations as well (see \cite{NY3}). 
One important point of symmetric forms is 
that the structure of B\"acklund transformations of these Painlev\'e 
equations can be described in a unified manner, by introducing a class 
of representations of affine Weyl groups inside certain Cremona groups. 
Also, with the $\tau$-functions appropriately defined, 
the dependent variables of the Painlev\'e equations allow certain 
``multiplicative formulas'' in terms of $\tau$-functions. 
One remarkable fact about our multiplicative formulas (\ref{f-by-tau})
is that the  factors are completely 
determined by the Cartan matrix of the corresponding 
affine root system.
Similar structures can be found commonly in 
various (discrete) integrable systems with 
Painlev\'e (singularity confinement) property
(\cite{RGH},\cite{KNS},\cite{KNH}).

\medskip
The main purpose of this paper is to present a new class of 
representations of affine Weyl groups which provides 
a prototype of affine Weyl group symmetry in nonlinear 
differential and difference systems.  

In Sections 1 and 2, we introduce a class of representations of the 
Coxeter groups 
of Kac-Moody type on certain fields of rational functions
(on the levels of $f$-{\em variables} and $\tau$-{\em functions}, 
respectively). 
This class of representations was found as a generalization of 
the structure of B\"acklund transformations in the symmetric forms 
of Painlev\'e equations $P_{\,\text{IV}}$, $P_{\,\text{V}}$ and 
$P_{\,\text{VI}}$ which are the cases of $A^{(1)}_2$, $A^{(1)}_3$ and 
$D^{(1)}_4$ respectively.

Our representation in the case of an affine root system 
provides naturally a discrete dynamical system from the lattice part 
of the affine Weyl group.  
We introduce in Section 3 the discrete dynamical 
systems associated with affine root systems in this sense. 
The case of $A^{(1)}_l$ is discussed in Section 4 in some detail as an example.
One interesting aspect of our system is that 
{\em continued fractions} arise naturally in the discrete dynamical system,  
with variations depending on the affine root system.  

In the final section, we explain how one can apply our discrete 
dynamical systems to the problem of symmetry of 
nonlinear differential (or difference) systems.  
In particular, we present a series of nonlinear ordinary differential systems
which have symmetry under the affine Weyl groups of type $A^{(1)}_l$.   
This series of nonlinear equations gives a generalization of 
the  Painlev\'e equations $P_{\,\text{IV}}$ and $P_{\,\text{V}}$ to higher 
orders. 

\section{A representation of the Coxeter group $W(A)$}

We fix a {\it generalized Cartan matrix} (or a {\it root datum}) 
$A=(a_{ij})_{i,j\in I}$ with $I$ being a finite indexing set. 
By definition, $A$ is a square matrix with the properties
\par\smallskip
\begin{quote}
(C1) \quad $a_{jj}=2 $ for all $j\in I$,\newline
(C2) \quad $a_{ij}$ is a nonpositive integer if $i\ne j$,\newline
(C3) \quad $a_{ij}=0 \Leftrightarrow a_{ji}=0$ \quad $(i,j\in I)$. 
\end{quote}
\par\smallskip\noindent
(See Kac \cite{Kac} for the basic properties of generalized Cartan matrices. 
Although we assume that $I$ is finite, 
a considerable part of the following argument can be formulated under the 
assumption that $A$ is {\em locally finite}, namely, for each $j\in I$, 
$a_{ij}=0 $  except for a finite number of $i$'s. )
We define the {\em root lattice} $Q=Q(A)$ 
and the {\em coroot lattice} $Q^\vee$ for $A$ by 
\begin{equation}
Q=\bigoplus_{j\in I} \BZ \, \alpha_j\quad\text{and}\quad
Q^\vee=\bigoplus_{j\in I} \BZ\, \alpha^\vee_j
\end{equation}
respectively, together with the pairing 
$\br{\,,\,} : Q^\vee\times Q \to \BZ$  
such that $\br{\alpha^\vee_i,\alpha_j}=a_{ij}$ for $i,j\in I$.
We denote by $W=W(A)$ 
the Coxeter group defined by the generators 
$s_i$ ($i\in I$) and defining relations 
\begin{equation}
s_i^2=1,\quad (s_is_j)^{m_{ij}}=1\quad(i,j\in I, i\ne j),
\end{equation}
where $m_{ij}=2,3,4,6$ or $\infty$ according as 
$a_{ij}a_{ji}=0,1,2,3$ or $\ge 4$. 
The generators $s_i$ act naturally on $Q$ by reflections 
\begin{equation}\label{s-on-a}
s_{i}(\alpha_j)=\alpha_j-\alpha_i\br{\alpha^\vee_i,\alpha_j}
=\alpha_j-\alpha_i a_{ij}
\end{equation}
for $i,j\in I$. 
Note that the action of each $s_i$ on $Q$ induces an 
automorphism of the field $\BC(\alpha)=\BC(\alpha_i; i\in I )$ 
of rational functions  in $\alpha_i$ $(i\in I)$ so that 
$\BC(\alpha)$ becomes a left $W$-module.

\medskip
Introducing a set of new ``variables'' $f_j$ ($j\in I$), we propose 
to extend the representation of $W$ on $\BC(\alpha)$ to the 
field $\BC(\alpha;f)=\BC(\alpha)(f_j; j\in I)$ of rational 
functions in $\alpha_j$ and $f_j$ ($j\in I$). 
In order to specify the action of $s_i$ on $f_j$, we fix 
a matrix $U=(u_{ij})_{i,j\in I}$ with entries in $\BC$ such 
that
\par\smallskip
\begin{tabular}{clll}
(0) &$u_{ij}=0$ & if & $i=j$ \ or \  $a_{ij}=0$, \\
(1) &$u_{ij}=-u_{ji}$ &if & $(a_{ij},a_{ji})=(-1,-1)$, \\
(2) &$u_{ij}=-u_{ji}$ or $-2u_{ji}$ & if & $(a_{ij},a_{ji})=(-2,-1)$, \\
(3) &$u_{ij}=-u_{ji},-\frac{3}{2}u_{ji},-2u_{ji}$ or $-3u_{ji}$ & if & 
$(a_{ij},a_{ji})=(-3,-1)$. 
\end{tabular}
\par\smallskip
\begin{thm}\label{thmA}
Let $A=(a_{ij})_{i,j\in I}$ be a generalized Cartan matrix and 
$U=(u_{ij})_{i,j\in I}$ a matrix satisfying the conditions above. 
For each $i\in I$, we extend the action of $s_i$ on $\BC(\alpha)$ to an 
automorphism of $\BC(\alpha;f)$ such that
\begin{equation}
s_i(f_j)=f_j+\frac{\alpha_i}{f_i}u_{ij} \quad(j\in I). 
\end{equation}
Then the actions of these $s_i$ define a representation of 
the Coxeter group $W=W(A)$ $($i.e. a left $W$-module structure$)$ 
on the field $\BC(\alpha;f)$ of rational functions.
\end{thm}
\noindent
We have only to check that the automorphisms $s_i$ on $\BC(\alpha;f)$
are involutions ($s_i^2=1$ for all $i\in I$) and that they satisfy 
the Coxeter relations $(s_is_j)^{m_{ij}}=1$ when $i\ne j$ and $m_{ij}=2,3,4,6$.
This can be carried out by direct computations since, 
for any $i\in I$, the automorphism $s_i$ stabilizes the subfield 
$\BC(\alpha)(f_i,f_k)$ for each $k\in I$ and, 
for any $i,j\in I$, both $s_i$ and $s_j$ stabilize 
the subfield $\BC(\alpha)(f_i,f_j,f_k)$ for each $k\in I$. 

\medskip
We remark that Theorem \ref{thmA} provides a systematic method to 
realize the Coxeter groups of Kac-Moody type nontrivially 
inside {\em Cremona groups} 
(groups of the birational transformations of affine spaces).

\begin{rem}
An important class of generalized Cartan matrices is that of 
{symmetrizable} ones, which includes the matrices of 
finite type and of affine type. 
Our condition on $U=(u_{ij})_{ij\in I}$
described above requires that $U$ should be ``almost'' skew-symmetrizable. 
The matrix $U$ can be thought of as specifying a sort of {\em orientation} of 
the Coxeter graph of $A$. 
It is also related to {\em Poisson structures} of dynamical systems. 
\end{rem}

\begin{rem}
Practically, it is sometimes necessary to consider the extension 
$\widetilde{W}=W\rtimes \Omega$ of 
$W=W(A)$ by a group $\Omega$ of diagram automorphisms 
of $A$. 
Recall that a {\em diagram automorphism} $\omega$ is by definition 
a bijection on $I$ such that 
$a_{\omega(i) \omega(j)}=a_{ij}$ for all $i,j\in I$;  
the commutation relations of each $\omega\in \Omega$
with elements of $W$ are given by  
$\omega s_i=s_{\omega(i)} \omega$ for all $i\in I$. 
Suppose that the matrix $U$ satisfies in addition 
the following compatibility condition with respect to $\Omega$ :  
$u_{\omega(i)\omega(j)}=u_{ij}$ for all $i,j\in I$, $\omega\in\Omega.$ 
Then, together with the automorphisms $\omega$ of $\BC(\alpha;f)$
such that 
$\omega(\alpha_j)=\alpha_{\omega(j)}$,
$\omega(f_j)=f_{\omega(j)}$  ($j\in I$),
the representation of $W$ in Theorem \ref{thmA} lifts to 
a representation of the {\em extended } Coxeter group 
$\widetilde{W}=W\rtimes \Omega$ on $\BC(\alpha;f)$.
\end{rem}

\section{$\tau$-Functions -- A further extension of the representation}

We now introduce another set of variables $\tau_j$ ($j\in I$), which 
we call the ``$\tau$-functions'' for 
the $f$-variables $f_j$ ($j\in I$).  
Considering the field extension  
$\BC(\alpha;f;\tau)$
$=\BC(\alpha;f)(\tau_j; j\in I)$,
we propose a way to extend the representation of $W$ 
of Theorem \ref{thmA} to $\BC(\alpha;f;\tau)$.

\begin{thm}\label{thmB} 
Let $A$ be a generalized Cartan matrix and $U=(u_{ij})_{i,j\in I}$ 
a matrix with entries in $\BC$ satisfying the conditions 
\begin{quote}
$(0)$ \quad $u_{jj}=0$ \quad for all $j\in I$, \newline
$(1)$ \quad $u_{ij}=u_{ji}=0$\quad if \ \ $a_{ij}=a_{ji}=0$, \newline
$(2)$ \quad $u_{ij}=-k u_{ji}$ \quad if \ \ $(a_{ij},a_{ji})=(-k,-1)$ with $k=1,2$ or $3$. 
\end{quote}
We extend the action of each generator $s_i$ of $W$ on 
$\BC(\alpha;f)$ to an automorphism of $\BC(\alpha;f;\tau)$
by the formulas 
\begin{equation}\label{s-on-tau}
s_i(\tau_j)=\tau_j\quad (i\ne j),\quad
s_i(\tau_i)=f_i \ \tau_i\prod_{k\in I} \tau_k^{-a_{ki}}=
f_i\frac{\prod_{k\in I\backslash\{i\}} \tau_k^{|a_{ki}|}}{\tau_i},
\end{equation}
for all $i,j\in I$.  
Then these automorphisms define a representation of $W$ on 
$\BC(\alpha;f;\tau)$
\end{thm}
The formulas (\ref{s-on-tau}) of Theorem \ref{thmB} specify 
how the $f$-variables should be expressed in terms of the $\tau$-functions:
\begin{equation} \label{f-by-tau}
f_j=\frac{\tau_j \ s_j(\tau_j)}
{{\prod}_{ i\in I\backslash\{j\}}\, \tau_i^{|a_{ij}|}}
\end{equation}
for all $j\in I$.  
We remark that this type of {\em multiplicative formulas by $\tau$-functions} 
is of a {\em universal\/} nature as can be found 
in various discretized integrable systems such as 
$T$-systems, discrete Toda equations and discrete Painlev\'e equations 
(see \cite{KNS},\cite{KNH},\cite{RGH},$\ldots$). 
In that context, the existence of multiplicative formulas is thought of 
as a reflection
of {\em singularity confinement} which is a discrete analogue of the 
Painlev\'e property.  

\begin{rem}
If the matrix $U$ is invariant with respect to a group $\Omega$ of 
diagram automorphisms, then the action of the extended Coxeter group 
$\widetilde{W}=W\rtimes \Omega$ on $\BC(\alpha;f)$ extends naturally 
to $\BC(\alpha;f;\tau)$ by $\omega.\tau_j=\tau_{\omega(j)}$ 
for all $j\in I$.
\end{rem}

\par\medskip 
Theorem \ref{thmB} can be proved essentially by direct computation to  
verify the fundamental relations of the Coxeter group 
with respect to the action on the $\tau$-functions $\tau_k$ ($k\in I$).  
Instead of giving the detail of such a proof, 
we will explain some of the ideas behind these multiplicative formulas. 
We consider that the $\tau$-functions should correspond to the 
{\em fundamental weights} $\Lambda_j$, while the $f$-variables do to 
{\em simple roots} $\alpha_i$.  
Let us denote by $L=\operatorname{Hom}_\BZ(Q^\vee,\BZ)$ the dual 
$\BZ$-module of the coroot lattice $Q^\vee$, and take the 
dual basis $\{\Lambda_j\}_{j\in I}$ of $\{\alpha^\vee_i\}_{i\in I}$ 
so that $L=\bigoplus_{j\in I} \BZ \Lambda_j$. 
Note that $L$, being the dual of $Q$, has a natural action of $W$
and that there is a natural $W$-homomorphism $Q\to L$ such that
\begin{equation}\label{a-in-L}
\alpha_j\mapsto\sum_{i\in I} \Lambda_i a_{ij}\quad(j\in I)
\end{equation} 
through the pairing $\br{\ ,\ }$.
(The lattice $L$ is in fact the weight lattice modulo the null roots.) 
The action of $W$ on $L$ is then described as 
\begin{equation}\label{s-on-L} 
s_i(\Lambda_j)=0\ (i\ne j), \quad
s_i(\Lambda_i)=\Lambda_i-\sum_{k\in I}\Lambda_k a_{ki}
%=-\Lambda_i+\sum_{k\in I;\,k\ne i} \Lambda_k |a_{ki}|
\end{equation}
for $i,j\in I$.  
We remark that formulas (\ref{s-on-tau}) in Theorem \ref{thmB} are 
a multiplicative analogue of $(\ref{s-on-L})$ {\em except for the 
factor} $f_j$.

Let us introduce the notation of formal exponentials for $\tau$-functions:
\begin{equation}
\tau^\lambda=\prod_{i\in I} \tau_i^{\lambda_i}\quad\text{for each}\quad 
\lambda=\sum_{i\in I}\lambda_i \Lambda_i\in L,
\end{equation}
where $\lambda_i=\br{\alpha_i^\vee,\lambda}$.
In order to clarify the meaning of Theorem \ref{thmB}, 
we consider the action of each element 
$w\in W$ on $\tau^\lambda$ for $\lambda\in L$.
Suppose now that the action of $W$ on $\BC(\alpha;f)$ can be 
extended to $\BC(\alpha;f;\tau)$ as described in Theorem \ref{thmB}. 
Since formulas (\ref{s-on-tau}) read as
$s_i(\tau^{\Lambda_j})= f_j^{\delta_{ij}} \,\tau^{s_i(\Lambda_j)}$
for $j\in I$,
we have by linearity 
\begin{equation}
s_i(\tau^{\lambda})= 
f_i^{\lambda_i} \,\tau^{s_i(\lambda)}
\end{equation}
for each $\lambda\in L$.
Hence, for each $w\in W$, we should have  rational functions 
$\phi_w(\lambda)\in \BC(\alpha;f)$ indexed by $\lambda\in L$ such that 
\begin{equation}\label{coboundary}
w(\tau^\lambda)=\phi_w(\lambda)\, \tau^{w.\lambda}\quad(w\in W, \lambda\in L).
\end{equation}
Furthermore, these functions $\phi_w(\lambda)$ should satisfy the 
following {\em cocycle condition}:
\begin{equation}\label{cocycle}
\phi_{w_1w_2}(\lambda)=w_1(\phi_{w_2}(\lambda))\, \phi_{w_1}({w_2}.\lambda)
\end{equation} 
for all $w_1,w_2\in W$ and $\lambda\in L$. 
Conversely, if one has a family $(\phi_w(\lambda))_{w\in W, \lambda\in L}$
of rational functions
satisfying the cocycle condition (\ref{cocycle}), 
one can define a representation of $W$ on $\BC(\alpha;f;\tau)$ by means of 
(\ref{coboundary}). 
Theorem \ref{thmB} is thus equivalent to the following proposition.
\begin{prop}\label{prop-cocycle}
Under the same assumption of Theorem \ref{thmB}, there exists a
unique cocycle $\phi=(\phi_w(\lambda))_{w\in W, \lambda\in L}$ such that 
\begin{equation}
\phi_{1}(\lambda)=1,\quad \phi_{s_i}(\lambda)=f_i^{\br{\alpha_i^\vee,\lambda}}
\quad(\lambda\in L)
\end{equation}
for each $i\in I$.
\end{prop}

\begin{rem}
Any family $\{\phi_w(\lambda)\}_{w\in W, \lambda\in L}$ 
of rational functions in $\BC(\alpha;f)$ can be identified with 
a mapping 
\begin{equation}\label{cochain}
\phi: W\to \operatorname{Hom}_\BZ(L,\BC(\alpha;f)^\times): 
w \mapsto \phi_w,
\end{equation}
where $\BC(\alpha;f)^\times$ stands for the multiplicative group of 
$\BC(\alpha;f)$ regarded as a $\BZ$-module.
The cocycle condition (\ref{cocycle}) is then equivalent to 
saying that the mapping $\phi$ of (\ref{cochain}) is a 
{\em Hochschild 1-cocycle}
of $W$ with respect to the natural $W$-bimodule structure of
$\operatorname{Hom}_\BZ(L,\BC(\alpha;f)^\times)$. 
Furthermore, formula (\ref{coboundary}) means that, 
this cocycle $\phi$ becomes the coboundary of the 0-cochain  
\begin{equation}
\tau \in \operatorname{Hom}_\BZ(L,\BC(\alpha;f;\tau)^\times): 
\lambda\mapsto \tau^\lambda 
\end{equation}
after the extension of the $W$-module $\BC(\alpha;f)$ 
to $\BC(\alpha;f;\tau)$. 
Thus one could say that: {\em The role of $\tau$-functions is to 
trivialize the Hochschild 1-cocycle defined by the 
$f$-variables.} 
From the cocycle condition, it follows that 
the cocycle $\phi_w : L\to \BC(\alpha;f)^\times$ 
of Proposition \ref{prop-cocycle} 
can be expressed as 
\begin{equation}
\phi_w(\lambda)=\prod_{r=1}^p s_{j_1}\cdots s_{j_{r-1}}
(f_{j_r})^{\br{\alpha^\vee_{j_r},\,
s_{j_{r+1}}\ldots s_{j_p}.\lambda}}
\quad(\lambda\in L)
\end{equation}
for any expression $w=s_{j_1}\ldots s_{j_p}$ of $w$ in terms of generators. 
\end{rem}

\medskip
The cocycle $\phi=(\phi_w(\lambda))_{w\in W,\lambda\in L}$ 
defined above plays a crucial role in application of our 
representation to discrete dynamical systems.  
One remarkable thing about this cocycle  is that $\phi$ seem to have 
a very strong {\em regularity} as described in the following conjecture. 
\begin{conj}\label{conj-cocycle}
In addition to conditions $(1)$ and $(2)$ of Theorem \ref{thmB}, 
suppose that the matrix $U=(u_{ij})_{i,j\in I}$ satisfies the condition 
\begin{quote}
$(3')$ $\quad u_{ij}a_{ji}+a_{ij}u_{ji}=0$ \quad for all $i,j\in I$. 
\end{quote}
Then, for any $k\in I$, the rational functions $\phi_w(\Lambda_k)$ $(w\in W)$ 
of $(\ref{coboundary})$ are polynomials in 
$\alpha_j,f_j$ and $u_{ij}$ $(i,j\in I)$ with coefficients in $\BZ$.
\end{conj}

\begin{rem}
In Sections 1 and 2, we presented a nontrivial class of representations
of Coxeter groups $W(A)$ over the fields of $f$-variables and 
$\tau$-functions, with $A$ being a generalized Cartan matrix. 
This class of representations appears in fact 
as {\em B\"acklund transformations}
(or the {\em Schlesinger transformations}) of the 
Painlev\'e equations $P_{\,\text{IV}}, P_{\,\text{V}}$ and $P_{\,\text{VI}}$,
which correspond to the cases of the generalized Cartan matrices $A$ of 
type $A^{(1)}_2$, $A^{(1)}_3$ and $D^{(1)}_4$, respectively. 
As to these Painlev\'e equations, one can define appropriate 
$f$-variables and $\tau$-functions 
for which the B\"acklund transformations are described as in 
Theorems \ref{thmA} and \ref{thmB} (see \cite{NY1}, \cite{NY3}).
(In the cases of $P_{\,\text{II}}$ and $P_{\,\text{III}}$,
which have symmetries of type $A^{(1)}_1$ and $C^{(1)}_2$,  
the corresponding representations of the affine Weyl groups 
on $f$-variables must be modified appropriately, while the 
multiplicative formulas in terms of $\tau$-functions keep the same structure.)
In the context of B\"acklund transformations of Painlev\'e equations, 
the functions $\phi_w(\Lambda_k)$, specialized to certain particular 
solutions, give rise to the special polynomials, called {\em Umemura 
polynomials} (see \cite{U},\cite{NOOU}), which are defined 
to be the main factors of $\tau$-functions for algebraic solutions 
of the Painlev\'e equations. 
For this reason, we expect that the functions $\phi_w(\Lambda_k)$ 
($w\in W, k\in I$) should supply an ample generalization of Umemura 
polynomials in terms of root systems.  
\end{rem}

\section{Affine Weyl groups and discrete dynamical systems}

In what follows, we assume that the generalized Cartan matrix $A$ is 
indecomposable and is {\em of affine type}. 
We use the standard notation of 
the indexing set $I=\{0,1,\ldots,l\}$ so that $\alpha_1,\ldots,\alpha_l$ 
form a basis for the corresponding finite root system. 
Recall that the null root $\delta$ is expressed as 
$\delta=a_0\alpha_0+a_1\alpha_1+\cdots+a_l\alpha_l$ with 
certain positive integers $a_0,a_1,\ldots,a_l$.  
The affine Weyl group $W=W(A)$ is generated by the 
fundamental reflections $s_0,s_1,\ldots,s_l$ with respect to the 
simple roots $\alpha_0,\alpha_1,\ldots,\alpha_l$:
\begin{equation}
W=W(A)=\br{s_0,\ldots,s_l}.  
\end{equation}
One important aspect of the affine case is that the affine Weyl group 
$W=W(A)$ has an alternative description as the semi-direct product of 
a free $\BZ$-submodule $M$ of rank $l$ 
of ${\overset\circ{\frak h}}_{\BR}=\bigoplus_{i=1}^l \BR \alpha_i^\vee$, 
and 
the finite Weyl group $W_0$ acting on $M$: 
\begin{equation}\label{MxW0} 
W\ \overset{\sim}{\leftarrow}\ M\rtimes W_0 \quad\text{with}\quad W_0
=\br{s_1,\ldots,s_l}. 
\end{equation}
For each element $\mu\in M$, we denote by $t_{\mu}$ the corresponding 
element of the affine Weyl group $W$, so that $t_{\mu+\nu}=t_\mu t_\nu$ 
for all $\mu,\nu\in M$. 
Note that the structure of the {\em lattice part} $M$  depends on the 
type of the affine root system and that, if $A$ is nontwisted, i.e., of type 
$X^{(1)}_l$, then $M$ is identified with the coroot lattice 
$\overset{\circ}{Q}{}^\vee$
of the finite root system with basis $\{\alpha_1,\ldots,\alpha_l\}$.
(There are descriptions analogous to (\ref{MxW0}) 
for certain {\em extended} affine Weyl groups 
$\widetilde{W}=W\rtimes\Omega$ as well.)
As we already remarked, the field 
$\BC(\alpha)=\BC(\alpha_0,\alpha_1,\ldots,\alpha_l)$ 
has a natural structure of $W$-module.  
The lattice part $M$ in the decomposition (\ref{MxW0}) acts 
on $\BC(\alpha)$ by 
\begin{equation}
t_\mu(\alpha_j)=\alpha_j-\br{\mu,\alpha_j}\delta\quad(j=0,1,\ldots,l;\mu\in M)
\end{equation}
as shift operators with respect to the simple affine roots. 
(The null root $\delta$ is a $W$-invariant element of $\BC(\alpha)$.  
For this reason, it is sometimes more 
convenient to consider $\delta$ to be a nonzero constant 
which represents the scaling of the lattice $M$.) 

\medskip
Suppose now that one has extended the action of $W$ from 
$\BC(\alpha)$ to $\BC(\alpha;f)=\BC(\alpha)(f_0,f_1,\ldots,f_l)$. 
At this moment, we can consider an arbitrary extension $\BC(\alpha;f)$ 
as a $W$-module, assuming that each element of $W$ acts on the function 
field as an automorphism; the representation of $W$ 
presented in Sections 1 and 2 provides a choice of such an extension. 
For each $\nu\in M$, 
we define a family of rational functions  
$F_{\nu j}(\alpha;f)\in \BC(\alpha; f)$ by 
\begin{equation}\label{disc-time-ev}
t_\nu(f_j)=F_{\nu j}(\alpha;f)\quad (j=0,1,\ldots,l). 
\end{equation} 
Then these formulas can already be considered as a 
{\em discrete dynamical system},
defined by a set of commuting discrete time evolutions. 
In other words, we obtain a commuting family of rational mappings
on the affine space
where $\alpha_j$ and $f_j$ play the role of coordinates 
of the discrete time variables and the dependent variables, 
respectively. 

To make clear the meaning of (\ref{disc-time-ev}) as a difference system, 
we set 
\begin{equation}
\alpha_j[\mu]=t_\mu(\alpha_j)=\alpha_j-\br{\mu,\alpha_j}\delta, 
\quad f_j[\mu]=t_\mu(f_j)\quad 
(j=0,\ldots,l)
\end{equation}
for each $\mu\in M$, 
and consider them as representing functions on $M$ with initial values 
$\alpha_j[0]=\alpha_j, f_j[0]=f_j$ ($j=0,\ldots,l$).  
Then formulas (\ref{disc-time-ev}) implies that
\begin{equation}\label{disc-time-ev2}
f_j[\mu+\nu]=
F_{\nu j}(\alpha[\mu];f[\mu])\quad(j=0,1,\ldots,l). 
\end{equation}
In this sense, the functions $F_{\nu j}(\alpha;f)$ defined above provide a 
difference dynamical system on the lattice $M$. 
Since  $f_j[\mu]$ is a rational function 
in $f_0,\ldots,f_l$, for each $\mu\in M$, 
the general solution of the difference system 
(\ref{disc-time-ev2})  {\em a priori\/} depends rationally
on initial values $f_0,f_1,\ldots,f_l$. 
Note also that the action of the affine Weyl group $W$ on $f_j[\nu]$ 
is described as 
\begin{equation}
(w.f_j)[\mu]=w(f_j[w^{-1}\mu])\quad(j=0,\ldots,l; \mu\in M)
\end{equation}
for all $w\in W$.  
In this sense, our difference system admits the action 
of the affine Weyl group $W(A)$. 
Note that, if one take the representation of Theorem \ref{thmA}, 
one has
\begin{equation}
(s_i.f_j)[\mu]
=f_j[\mu]+\frac{\alpha_i[\mu]}{f_i[\mu]}u_{ij}
\end{equation}
for $i,j=0,\ldots,l$. 

Suppose that one can extend the action of $W$ further to the $\tau$-functions 
as in Theorem \ref{thmB} and set 
$\tau_i[\mu]=t_\mu(\tau_i)$, regarding $\tau_i[0]=\tau_i$ as initial 
values of the $\tau$-functions. 
Then from (\ref{f-by-tau}) we obtain the multiplicative formulas 
\begin{equation} 
f_j[\mu]=\frac{\tau_j[\mu] \ s_{\alpha_j[\mu]}(\tau_j[\mu])}
{{\prod}_{ i\in I\backslash\{j\} }\, \tau_i[\mu]^{|a_{ij}|}}
\qquad(j=0,\ldots,l)
\end{equation}
for the $f$-variables in terms of $\tau$-functions.  
In terms of the cocycle $\phi$, these formulas are rewritten 
by  (\ref{coboundary}) into 
\begin{equation}
f_j[\mu]=\frac{\phi_{t_\mu}(\Lambda_j)\ \phi_{t_\mu s_j}(\Lambda_j)}
{{\prod}_{ i\in I\backslash\{j\}}\, \phi_{t_\mu}(\Lambda_i)^{|a_{ij}|}}
\qquad(j=0,\ldots,l), 
\end{equation}
which give a complete description of the general solution of the 
difference system (\ref{disc-time-ev}) in terms of the initial 
values $f_0,\ldots,f_l$.  
In this sense,  the cocycle $\phi$ {\em solves} 
our difference system (\ref{disc-time-ev}).  
It should be noted that all these properties of the difference system 
(\ref{disc-time-ev}), or  (\ref{disc-time-ev2}) equivalently, 
are already guaranteed when we take the representation of the affine 
Weyl group $W(A)$ as in Theorem \ref{thmB}.  
Also, it is meaningful if one could find other types of representations 
of affine Weyl groups which have the properties of Theorem \ref{thmB}.  

\medskip
We now take the representation of the affine Weyl group 
$W(A)$ on $\BC(\alpha;f;\tau)$ introduced in Theorem \ref{thmB}. 
One interesting feature of our representation is that 
{\em continued fractions} arise naturally in the description 
of discrete dynamical systems, and that the structure of continued 
fractions is determined by the affine root system. 
We assume for simplicity that the generalized Cartan matrix $A$ is 
of type $X^{(1)}_l$.  
For a given element $w\in W(A)$, take a reduced decomposition 
$w=s_{i_1}s_{i_2}\cdots s_{i_p}$ of $w$, and define the affine 
roots $\beta_1,\beta_2,\cdots,\beta_p$ by 
\begin{equation}
\beta_1=\alpha_1, \ \ \beta_2=s_{i_1}(\alpha_{i_2}), \ \ 
\ldots, \ \ \beta_p=s_{i_1}\cdots s_{i_{p-1}}(\alpha_{i_p}).  
\end{equation}
Note that these $\beta_r$ ($r=1,\ldots,p$)  give precisely the set 
of all positive real roots whose reflection hyperplanes separate 
the fundamental alcove $C$ and its image $w.C$ by $w$.  
Since the action of $s_i$ on $f_j$ is given by 
\begin{equation}
s_i(f_j)=f_j+\frac{\alpha_i}{f_i} u_{ij}\qquad(i,j=0,\ldots,l),
\end{equation}
we have inductively
\begin{equation}\label{w-on-f}
w(f_j)= f_j + \frac{\alpha_{i_1}}{f_{i_1}}u_{i_1j}
+s_{i_1}\big(\frac{\alpha_{i_2}}{f_{i_2}}\big)u_{i_2j}
+\cdots+s_{i_1}\cdots s_{i_{p-1}}
\big(\frac{\alpha_{i_p}}{f_{i_p}}\big)u_{i_pj}. 
\end{equation}
Each summand of this expression is given by the continued 
fraction 
\begin{equation}
s_{i_1}\cdots s_{i_{r-1}}\big(\frac{\alpha_{i_r}}{f_{i_r}}\big)
=\cfrac{\beta_r}{f_{i_r}+
u_{i_{r-1}i_{r}}
\cfrac{\beta_{r-1}}{
\overset{\displaystyle{f_{i_{r-1}}+{}_{\ddots}}\phantom{\frac{X}{X}}\hfill}
{\phantom{SSSSSS}\cfrac{\beta_2}{
{f_{i_2}+u_{i_1i_2}\dfrac{\beta_{1}}{ f_{i_1}}}
}}}}
\end{equation}
along the reduced decomposition $w=s_{i_1}\cdots s_{i_p}$. 
Note also that formula (\ref{w-on-f}) for $w(f_j)$ has an alternative 
expression 
\begin{equation}
w(f_j)=\frac{\phi_{w}(\Lambda_j)\ \phi_{ws_j}(\Lambda_j)}
{\prod_{i\in I\backslash\{j\}} \phi_w(\Lambda_i)^{|a_{ij}|}}
\end{equation}
in terms of the cocycle $\phi$, which is implied by 
(\ref{f-by-tau}) and (\ref{coboundary}).  
If one take an element $\nu\in M=\overset{\circ}{Q}{}^\vee$ of the dual root lattice,  
the rational functions $t_{\nu}(f_j)=F_{\nu j}(\alpha;f)$ 
$(j=0,\ldots,l)$ for the 
time evolution with respect to $\nu$ are 
determined in the form 
\begin{equation}
F_{\nu j}(\alpha;f)=f_j+\sum_{r=1}^p 
s_{i_1}\cdots s_{i_{r-1}}\big(\frac{\alpha_{i_r}}{f_{i_r}}\big) u_{i_rj}
\end{equation}
as a sum of continued fractions 
along the reduced decomposition of $t_{\nu}$, 
with positive real roots separating the fundamental alcove $C$ and 
its translation $C+\nu$.
We remark that a similar description of the rational functions 
$F_{\nu j}(\alpha;f)$  can be given also for the cases of 
extended affine Weyl groups $\widetilde{W}=W\rtimes \Omega$. 
A series of such discrete dynamical systems 
will be given in the next section. 

\section{Discrete dynamical system of type $A^{(1)}_l$}
As an example of our discrete dynamical systems associated with 
affine root systems, we will give an explicit description of the 
case of $A^{(1)}_l$ with $l\ge 2$. 
Consider the generalized Cartan matrix 
\begin{equation}
A=\begin{pmatrix}
2 & -1 &  0  &\cdots &0 & -1 \\
-1& 2  & -1 &\cdots &0 &  0 \\
0 & -1  & 2 &\cdots &0 &  0 \\
\vdots & \vdots  & \vdots &\ddots & \vdots & \vdots \\
0 & 0 & 0 & \cdots & 2 &-1 \\
-1& 0  & 0  &\cdots  & -1 & 2
\end{pmatrix}
\end{equation}
of type $A^{(1)}_l$ ($l\ge 2$), and 
identify the indexing set $\{0,1,\ldots,l\}$ with $\BZ/(l+1)\BZ$. 
We take following matrix of ``orientation'' to specify our 
representation of $W=W(A^{(1)}_l)$: 
\begin{equation}\label{AU}
U=\begin{pmatrix}
0 & 1 &  0  &\cdots &0 & -1 \\
-1& 0  & 1 &\cdots &0 &  0 \\
0 & -1  & 0 &\cdots &0 &  0 \\
\vdots & \vdots  & \vdots &\ddots & \vdots & \vdots \\
0 & 0 & 0 & \cdots & 0 & 1 \\
1& 0  & 0  &\cdots  & -1 & 0
\end{pmatrix}.  
\end{equation}
Then the action of the affine Weyl group $W=\br{s_0,\ldots,s_l}$ on the 
variables $\alpha_j$, $f_j$ and $\tau_j$ is given explicitly as follows:
\begin{equation}\label{WAl1}
\begin{array}{lllll}
s_i(\alpha_i)=-\alpha_i, 
&s_i(\alpha_{j})=\alpha_j+\alpha_i&(j=i\pm 1),
&s_i(\alpha_{j})=\alpha_j&(j\ne i,i\pm1 ) \\
s_i(f_i)=f_i,
&s_i(f_j)=\displaystyle{f_j\pm\frac{\alpha_i}{f_i}}&(j=i\pm 1),
&s_i(f_j)=f_j&(j\ne i,i\pm1) \\
s_i(\tau_i)=\displaystyle{f_i \,\frac{\tau_{i-1}\tau_{i+1}}{\tau_i}},
&s_i(\tau_j)=\tau_j &(j\ne i)
\end{array}
\end{equation}
Note that $U$ is invariant with respect to the diagram rotation 
$\pi: i\to i+1$. 
Hence this action of $W$ extends to 
the extended affine Weyl group
$\widetilde{W}=W\rtimes\{1,\pi,\ldots,\pi^l\}$ by 
\begin{equation}\label{WAl2}
\pi(\alpha_j)=\alpha_{j+1},\quad \pi(f_j)=f_{j+1},\quad \pi(\tau_j)=\tau_{j+1}.
\end{equation}
The group $\widetilde{W}$ is now 
isomorphic to $\overset{\circ}{P}\rtimes W_0$, 
where $\overset{\circ}{P}$ is the weight lattice of the finite root 
system of type $A_l$ and $W_0=\br{s_1,\ldots,s_l}\simeq {\frak S}_{l+1}$.
Taking the first fundamental weight 
$\varpi_1=(l\alpha_1+(l-1)\alpha_2+\cdots+\alpha_l)/(l+1)$ 
of the finite root system, we set
\begin{equation}
T_1=t_{\varpi_1},\quad T_i=\pi T_{i-1}\pi^{-1}\quad (i=2,\ldots,{l+1}).
\end{equation} 
These {\em shift operators\/} are expressed as 
\begin{equation}
T_1=\pi s_ls_{l-1}\cdots s_1, \ \ 
T_2=s_1\pi s_l\ldots s_2,\ \ \ldots,\ \  T_{l+1}=s_l\cdots s_1 \pi
\end{equation}
in terms of the generators of $\widetilde{W}$.
Note that $T_1\cdots T_{l+1}=1$ and that 
$T_1,\ldots,T_l$ form 
a basis for the lattice part of $\widetilde{W}$. 

The simple affine roots $\alpha_0,\ldots,\alpha_l$ are the dynamical variables 
for the shift operators $T_1,\ldots,T_{l}$ such that  
\begin{equation}
T_i(\alpha_{i-1})=\alpha_{i-1}+\delta,\quad
T_i(\alpha_{i})=\alpha_i-\delta, \quad
T_i(\alpha_j)=\alpha_j \ \  (j\ne i-1,i). 
\end{equation}
For each $k\in \BZ/(l+1)\BZ$ and $r=0,1,\ldots,l-1$, we define 
$g_{k,r}$ to be the continued fraction
\begin{equation}
\begin{align}
g_{k,r}&=s_{k+r}s_{k+r-1}\cdots s_{k+1}\big(\frac{\alpha_k}{f_k}\big)\\
&=\contfrac{\alpha_k+\ldots+\alpha_{k+r}}{f_k}-
\contfrac{\alpha_{k+1}+\ldots+\alpha_{k+r}}{f_{k+1}}-
\cdots-\contfrac{\alpha_{k+r}}{f_{k+r}}.
\notag
\end{align}
\end{equation}
Then the discrete time evolution by $T_1$ is expressed as
\begin{equation}\label{dAl}
\begin{align}
T_1(f_0)&=f_1-g_{2,l-1}+g_{0,0},\\
T_1(f_1)&=f_2-g_{3,l-2},\notag\\
T_1(f_2)&=f_3-g_{4,l-3}+g_{2,l-1},\notag\\
& \cdots\notag\\
T_1(f_{l-1})&=f_l-g_{0,0}+g_{l-1,2},\notag\\
T_1(f_l)&=f_0+g_{l,1}.\notag
\end{align} 
\end{equation}
The corresponding formulas for $T_2,\ldots, T_l$ are obtained 
from these by applying the diagram rotation $\pi$. 
\section{Nonlinear systems with affine Weyl group symmetry}

As we already remarked in Section 2, 
the representation of $W(A)$ introduced in Theorems \ref{thmA} and \ref{thmB}, 
for the cases $A^{(1)}_2$, $A^{(1)}_3$ and $D^{(1)}_4$, 
arises in nature as B\"acklund transformations of Painlev\'e equations 
$P_{\,\text{IV}}, P_{\,\text{V}}$ and $P_{\,\text{VI}}$, 
respectively.  
Hence, the Painlev\'e equations 
$P_{\,\text{IV}}, P_{\,\text{V}}$ and $P_{\,\text{VI}}$
have the structure of discrete dynamical systems on the lattice 
as described in Section 3 with respect to B\"acklund transformations. 
As to $P_{IV}$, this point has been discussed in detail in our 
previous paper \cite{NY1}. 
``Symmetric forms'' of all Painlev\'e equations 
$P_{\text{II}},\ldots, P_{\,\text{VI}}$ and their B\"acklund transformations 
will be discussed in our forthcoming paper \cite{NY3}. 

\medskip
From the viewpoint of nonlinear equations of Painlev\'e type, 
an important problem would be the following: 

\begin{prob}\label{prob}
For each affine root system $($or for each generalized Cartan matrix $A$, 
in general$)$, find a system of differential $($or difference$)$equations 
for which the Coxeter group $W=W(A)$ acts as B\"acklund transformations. 
\end{prob}

\noindent
We believe that such differential (or difference) systems 
with affine Weyl group symmetry 
should provide an intriguing class of dynamical systems with 
rich mathematical structures, to be compared to Painlev\'e equations. 
We also remark that, if one specifies the representation of $W=W(A)$ 
in advance as in Theorem \ref{thmB}, then 
the problem mentioned above is equivalent to finding such 
derivations (or shift operators) on $\BC(\alpha;f)$ and 
$\BC(\alpha;f;\tau)$ 
that commute with the action of $W(A)$. 

In this section, we will introduce some examples of type $A^{(1)}_l$ 
of difference and differential systems with affine Weyl group 
symmetry, as well as remarks on the continuum limit 
from the difference to the differential systems. 

\medskip
We first explain a general idea to construct 
difference systems with affine Weyl group symmetry
by means of our discrete dynamical systems associated with 
affine root systems.
Consider the 
discrete dynamical system defined by an affine root system 
as in Section 3.  
If we take a sublattice $N\subset M$ of rank $r$,  
then the centralizer $Z_{W(A)}(N)$ of $N$ in $W$ gives rise to a 
group of B\"acklund transformations of the discrete system
\begin{equation}\label{disc-time-ev-r} 
t_\nu(f_j)=F_{\nu j}(\alpha;f)\quad(j=0,\ldots,l)
\end{equation}
on the sublattice $N$ of rank $r$, with $\alpha_j$ ($j=0,\ldots,l)$ 
regarded as functions on $N$ such that 
$t_\nu(\alpha_j)=\alpha_j-\br{\nu,\alpha_j}$. 
The centralizer $Z_{W(A)}(N)$ contains in fact 
subgroups generated by reflections acting on the quotient $M/N$. 
For instance, let $W_{M/N}$ be the group generated by 
the reflections $s_\alpha$ with respect to the affine roots $\alpha$ 
that are perpendicular to the lattice $N$. 
Then $W_{M/N}$ is contained in the group of B\"acklund 
transformations of the discrete system (\ref{disc-time-ev-r}).
(The group of symmetry thus obtained 
may have a different structure from that of our 
representations of Sections 1 and 2.) 

For example, the difference system (\ref{dAl}) with respect to
the shift operator $T_1$ has symmetry under the affine Weyl 
group $W(A^{(1)}_{l-1})=\br{r,s_2,\ldots,s_l}$, where 
$r=s_0 s_1 s_0$.  
The corresponding simple affine roots are given by 
$\alpha_0+\alpha_1, \alpha_2,\ldots,\alpha_l$. 
Note that the root $\alpha_0+\alpha_1$ is invariant under $T_1$.
The reflection $r$ acts on the variables $f_j$ as follows:
\begin{equation}
\begin{align}
&r(f_1)=f_1+\frac{\alpha_0+\alpha_1}{s_1(f_0)},
\quad
r(f_2)=f_2+\frac{\alpha_0+\alpha_1}{s_0(f_1)},
\\
&r(f_l)=f_l-\frac{\alpha_0+\alpha_1}{s_1(f_0)},
\quad
r(f_0)=f_0-\frac{\alpha_0+\alpha_1}{s_0(f_1)}, \notag
\\
&r(f_j)=f_j \qquad(j=3,\ldots,l-1).\notag
\end{align}
\end{equation}
We remark 
that the two element $s_1(f_0)$ and $s_0(f_1)= T_1 s_1(f_0)$ 
are invariant under the action of $r$ as well as 
$f_3,\ldots,f_{l-1}$. 
Similarly, the difference system with respect to the commuting 
operators $T_1,\ldots, T_k$ has affine Weyl group symmetry under the 
subgroup $W(A^{(1)}_{l-k})=\br{r,s_{k+1},\ldots,s_{l}}$ where 
$r=s_0s_1\cdots s_{k-1} s_k s_{k-1}\cdots s_0$. 

\medskip
If one can take an appropriate continuum limit of 
the sublattice $N$ inside $M$, one would possibly obtain 
a differential system in $r$ variables 
whose group of B\"acklund transformations contains 
a reasonable reflection group. 
We show an example in which the idea explained above works nicely, 
in some detail. 

Consider the discrete dynamical system of type $A^{(1)}_2$ 
with extended affine Weyl group $\widetilde{W}=W\rtimes \Omega$
as in Section 4.  
We take the element $T=T_1=\pi s_2 s_1\in \widetilde{W}$ which represents 
to the translation $t_{\varpi_1}$ with respect to the first fundamental 
weight $\varpi_1=(2\alpha_1+\alpha_2)/3$ of the finite root system, so that 
\begin{equation}
T(\alpha_0)=\alpha_0+\delta,\quad T(\alpha_1)=\alpha_1-\delta,
\quad T(\alpha_2)=\alpha_2. 
\end{equation}
Note that $T$ can be considered as a shift operator with respect to 
the variable $\alpha_1$. 
Our discrete dynamical system for this case is described as follows:
\begin{equation}
\begin{align}\label{Tf0}
T(f_0)&=f_1+\frac{\alpha_0}{f_0}-
\contfrac{\alpha_0+\alpha_1}{f_2}-\contfrac{\alpha_0}{f_0},\\
T(f_1)&=f_2-\frac{\alpha_0}{f_0},\notag\\
T(f_2)&=f_0+\contfrac{\alpha_0+\alpha_1}{f_2}-\contfrac{\alpha_0}{f_0}\notag
\end{align}
\end{equation}
in terms of continued fractions.
We remark that $T^{-1}(f_0)$ takes a simpler form than $T(f_0)$ above:
\begin{equation}\label{Tf1}
T^{-1}(f_0)=f_2+\frac{\alpha_1}{f_1}.%=\frac{f_1f_2+\alpha_1}{f_1}.
\end{equation}
Noticing that the element $f_0+f_1+f_2$ is invariant under the 
action of $\widetilde{W}$, we set $f_0+f_1+f_2=c$.  
Then from (\ref{Tf0}) and (\ref{Tf1}) 
we  obtain the following equivalent form of our difference system: 
\begin{equation}\label{dP20}
T^{-1}(f_0)+f_0=c-f_1+\frac{\alpha_1}{f_1},\quad 
f_1+T(f_1)=c-f_0-\frac{\alpha_0}{f_0}. 
\end{equation}
With the notation $f_i[n]=T^n(f_i)$ for $n\in\BZ$, this equation gives rise 
to a representation of the {\em second discrete Painlev\'e equation\/}
$dP_{\,\text{II}}$ (\cite{SKGHR},\cite{S}): 
\begin{equation} 
\begin{align}\label{dP2}
f_0[n-1]+f_0[n]&=c-f_1[n]+\frac{\alpha_1-n\delta}{f_1[n]},\quad \\
f_1[n]+f_1[n+1]&=c-f_0[n]-\frac{\alpha_0+n\delta}{f_0[n]}\quad(n\in \BZ).\notag
\end{align}
\end{equation}
Since the shift operator $T=\pi s_2 s_1$ commute with the two reflections 
$r_0=s_0s_1s_0$ and $r_1=s_2$,  we see that the difference system 
(\ref{dP20})
or (\ref{dP2}) 
has symmetry of the affine Weyl group $W(A^{(1)}_1)=\br{r_0,r_1}$. 
(The corresponding simple roots are 
$\beta_0=\alpha_0+\alpha_1$ and $\beta_1=\alpha_2$.)

The second Painlev\'e equation $P_{\,\text{II}}$ 
arises as a continuum limit of the difference system (\ref{dP2}),
and that $A^{(1)}_1$-symmetry of (\ref{dP2}) naturally passes to 
$P_{\,\text{II}}$.
Introduce a small parameter $\varepsilon$ such that $\delta=\varepsilon^3$, 
and set 
\begin{equation}
\begin{align}
&f_0[n]=1+\varepsilon \psi + \varepsilon^2 \varphi_0,\quad
f_1[n]=1-\varepsilon \psi + \varepsilon^2 \varphi_1,\quad
c=2,\\
&\alpha_0+n \delta=-1+\varepsilon^2 x+\varepsilon^3 a_0,\ \ 
\alpha_1-n \delta=1-\varepsilon^2 x+\varepsilon^3 a_1,\ \ 
\alpha_2=\varepsilon^3 b_1. \notag
\end{align}
\end{equation} 
Then in the limit as $\varepsilon\to 0$, the difference equations (\ref{dP2}) 
imply the following differential equation for $\varphi_0, \varphi_1,\psi$: 
\begin{equation}
\begin{align}
&\varphi_0' = 2\varphi_0 \psi+a_0-\frac{1}{2},\  
\varphi_1' = 2\varphi_1 \psi+a_1-\frac{1}{2},\\
&\psi'=2(\varphi_0+\varphi_1)-\psi^2+x.\notag
\end{align}
\end{equation}
{}From this we get the second Painlev\'e equation for $\psi$ 
\begin{equation}
\psi''=2 \psi^3 - 2 x \psi -2 b_1+1,
\end{equation}
and the other dependent variables $\varphi_0, \varphi_1$ are 
determined by quadrature from $\psi$. 
At the same time, we obtain 
the following B\"acklund transformations $r_0$ and $r_1$ for $\psi$:
\begin{equation}
r_0(\psi)=\psi-\frac{2 b_0}{\psi'-\psi^2+x},\quad 
r_1(\psi)=\psi-\frac{2 b_1}{\psi'+\psi^2-x},
\end{equation}
where $b_0=a_0+a_1=1-b_1$.  
The parameters $b_0, b_1$ are the simple roots for the $A^{(1)}_{1}$-symmetry
of $P_{\,\text{II}}$. 

\medskip
Finally, 
we present a series of differential systems with $A^{(1)}_l$-symmetry 
$(l\ge 2)$, which give a generalization of the Painlev\'e equations 
$P_{\,\text{IV}}$ and $P_{\,\text{V}}$. 

In our previous paper \cite{NY1}, we introduced the symmetric form of the 
fourth Painlev\'e equation: 
\begin{equation}\label{SP4}
\begin{align}
&f'_0=f_0(f_1-f_2)+\alpha_0, \\
&f'_1=f_1(f_2-f_0)+\alpha_1, \notag\\
&f'_2=f_2(f_0-f_1)+\alpha_2. \notag 
\end{align}
\end{equation}
This system defines in fact a derivation $'$ of the field $\BC(\alpha;f)$ 
which commute with the action of the extended affine Weyl group
$\widetilde{W}$ of type $A^{(1)}_2$ as in (\ref{WAl1}) and (\ref{WAl2}).
(Note that the convention of \cite{NY1} corresponds to the 
transposition of $U$ in (\ref{AU}).) 
We remark that the sum $f_0+f_1+f_2$ is invariant under $\widetilde{W}$, 
and satisfies the equation $(f_0+f_1+f_2)'=\alpha_0+\alpha_1+\alpha_2=\delta$. 
Introduce the independent variable $x$ so that $x'=1$, and eliminate one of 
the three $f$-variables, noting that $f_0+f_1+f_2$ is a linear function of $x$.  
Then the differential system above is rewritten into a system of order 2, 
which is equivalent to the Painlev\'e equation $P_{\,\text{IV}}$. 

Differential system (\ref{SP4}) has a generalization to higher orders.  
For example, when $l=4$, the differential system 
\begin{equation}
\begin{align}
&f'_0=f_0(f_1-f_2+f_3-f_4)+\alpha_0, \\
&f'_1=f_1(f_2-f_3+f_4-f_0)+\alpha_1, \notag\\
&f'_2=f_2(f_3-f_4+f_0-f_1)+\alpha_2, \notag\\
&f'_3=f_3(f_4-f_0+f_1-f_2)+\alpha_3, \notag\\ 
&f'_4=f_4(f_0-f_1+f_2-f_3)+\alpha_4 \notag 
\end{align}
\end{equation}
has $A^{(1)}_4$-symmetry.  
Note that the sum $f_0+f_1+f_2+f_3+f_4$ is a linear function 
of the independent variable $x$ such that $x'=1$ 
and that the system above is essentially of order 4. 
In general, when $l=2n$, 
the following differential system (essentially of order $2n$) turns out to 
have $A^{(1)}_{2n}$-symmetry with the B\"acklund transformations 
defined as in Section 4: 
\begin{equation}
f_j'=f_j\sum_{1\le r\le n}(f_{j+2r-1}-f_{2r}\big)+\alpha_j
\quad (j=0,1,\ldots,2n).
\end{equation}
We remark that this differential system is obtained 
as a continuum limit from the difference system with 
$A^{(1)}_{2n}$-symmetry which arises from the discrete dynamical 
system of type $A^{(1)}_{2n+1}$, in the manner as we explained above. 

We also found a series of differential systems with $A^{(1)}_{2n+1}$-symmetry 
($n=1,2,\ldots$) which generalize the fifth Painlev\'e equation $P_{\,\text{V}}$: 
\begin{equation}\label{A2n+1}
\begin{align}
f_j'&=f_j\big(
\sum_{ 1\le r\le s\le n } f_{j+2r-1}f_{j+2s}
-\sum_{ 1\le r\le s\le n} f_{j+2r}f_{j+2s+1}
\big)\\
&\quad+(\frac{\delta}{2}-\sum_{1\le r\le n}\alpha_{j+2r})f_j +
\alpha_j(\sum_{1\le r\le n}f_{j+2r})
\quad (j=0,1,\ldots,2n+1),\notag
\end{align}
\end{equation}
where $\alpha_0+\cdots+\alpha_{2n+1}=\delta$. 
We remark that differential system (\ref{A2n+1}) is 
also essentially of order $2n$, 
since each of the sums $\sum_{r=0}^n f_{2r}$ and $\sum_{r=0}^n f_{2r+1}$ 
is determined elementarily. 
The Painlev\'e equation $P_{\,\text{V}}$ is covered as the case $n=1$
(see \cite{NY3}):
\begin{equation}
\begin{align}
f_0'&=f_0(f_1f_2-f_2f_3)+(\frac{\delta}{2}-\alpha_2)f_0+\alpha_0 f_2,
\\
f_1'&=f_1(f_2f_3-f_3f_0)+(\frac{\delta}{2}-\alpha_3)f_1+\alpha_1 f_3,
\notag\\
f_2'&=f_2(f_3f_0-f_0f_1)+(\frac{\delta}{2}-\alpha_0)f_2+\alpha_2 f_0,
\notag\\
f_3'&=f_3(f_0f_1-f_1f_2)+(\frac{\delta}{2}-\alpha_1)f_3+\alpha_3 f_1,
\notag
\end{align}
\end{equation}
where $\alpha_0+\alpha_1+\alpha_2+\alpha_3=\delta$. 

\medskip
These two series of differential systems with affine Weyl group symmetry 
can be considered as a variation of Lotka-Voltera equations and 
Bogoyavlensky lattices, including the 
parameters $\alpha_0,\ldots,\alpha_l$.  
Also, the structures of their B\"acklund transformations 
can be described completely in terms of the discrete dynamical 
systems we have introduced in this paper.
(Details will be discussed elsewhere.) 
We expect that these systems of differential equations  
with affine Weyl group symmetry 
deserve to be studied individually from various aspects, 
since they already give a candidate for 
systematic generalization of 
Painlev\'e equations to higher orders.  

\vfill\eject

\end{document}